\documentclass[12pt,a4paper]{article}
\usepackage{a4wide}
\usepackage{amsmath,amsfonts,amssymb,amsthm}
\usepackage{mathrsfs}
\usepackage{ifthen}
\usepackage{epsfig}

\advance\textwidth by 16mm
\advance\hoffset -8mm
\advance\textheight by 30mm
\advance\voffset -15mm

\newtheorem{Theorem}{Theorem}[section]
\newtheorem*{Conjecture}{Conjecture}

\newtheorem{Def}{Definition}[section]
\newtheorem{Lemma}{Lemma}[section]

\makeatletter
\@addtoreset{equation}{section}
\makeatother

\makeatletter
\newcommand\txt[2]{{#1{\ifthenelse{\equal {#2}{}}{}{ #2}}. }}
\def\Remark#1.#2\par{\par\medskip\noindent\txt{\bf Remark}{#1}\ #2\par%
               \@ifnextchar {\Remark}{}{\@ifnextchar\begin{}{\smallskip}}}
\def\Example#1.{\txt{\bf Example}{#1}}
\makeatother

\newcommand\emp\varnothing
\newcommand\eps\varepsilon

\def\^#1{^{\overline{#1}}}

\def\conv{\mathop{\operatorfont conv}}

\def\supp{\mathop{\operatorfont supp}\nolimits}

\def\pr{\mathop{\operatorfont pr}\nolimits}
\def\lchi{\mathop{\ell\chi}\nolimits}
\def\RR{{\mathbb R}}

\begin{document}

\title{Examples of topologically highly chromatic graphs with locally small chromatic number}

\author{Ilya I. Bogdanov\thanks{The work was supported by the Russian government project 11.G34.31.0053 and by RFBR grant No.~13-01-00563.}}

\maketitle

\begin{abstract}
  Kierstead, Szemer\'edi, and Trotter showed that a graph with at most $\lfloor r/(2n)\rfloor^n$ vertices such that each ball of radius $r$ in it is $c$-colorable should have chromatic number at most $n(c-1)+1$. We show that this estimate is sharp in $r$. Namely, for every $n$, $r$, and $c$ we construct a graph $G$ containing $O((2rc)^{n-1}c)$ vertices such that $\chi(G)\geq n(c-1)+1$, although each ball of radius $r$ in $G$ is $c$-colorable. The core idea is the construction of a graph whose neighborhood complex is homotopy equivalent to the join of neighborhood complexes of two given graphs.
\end{abstract}

\section{Introduction}

Let $G=(V,E)$ be a graph (with no loops or multiple edges). By $d_G(u,v)$ we denote the {\em distance} between the vertices $u,v\in V$. A subset $V_1\subseteq V$ is {\em independent} if none of the edges has both endpoints in $V_1$. The {\em chromatic number} $\chi(G)$ of~$G$ is the minimal number of colors in a proper coloring of~$G$, that is --- the minimal number of parts in a partition of $V$ into independent sets.

\begin{Def}
  Let $r$ be a positive integer. The {\em ball} of radius $r$ with center $v\in V$ is the set $U_r(v,G)=\{u\in G\colon d(u,v)\leq r\}$. The {\em $r$-local chromatic number} $\lchi_r(G)$ of a graph~$G$ is the maximal chromatic number of a ball of radius~$r$ in~$G$.
\end{Def}

Note that even for $r=1$ our definition of the local chromatic number is quite different from that introduced by Erd\H os et al. in~\cite{erd-loc}.

By the celebrated result of Erd\H os~\cite{erdos}, for every integer $n>1$ and $g>2$ there exists a graph of girth~$g$ and chromatic number greater than~$n$; thus for every $r$ there exist a graph~$G$ with $\lchi_r(G)=2$ and arbitrarily large $\chi(G)$. Later Erd\H os~\cite{erdos3} conjectured that for every positive integer~$s$ there exists a constant $c_s$ such that the chromatic number of each graph~$G$ having $N$ vertices and containing no odd cycles of length less than $c_sN^{1/s}$ does not exceed $s+1$. This conjecture was proved by Kierstead, Szemer\'edi, and Trotter~\cite{kst}. In fact, they have proved the following more general result.

\begin{Theorem}[{\cite[Theorem 1]{kst}}]
  \label{kst}
  Assume that $G=(V,E)$ is a graph such that $\lchi_r(G)\leq c$ and $|V|\leq \lfloor r/(2n)\rfloor^n$. Then $\chi(G)\leq n(c-1)+1$.
\end{Theorem}

They have also posed a question whether this bound is sharp. The strong form of this question is as follows.

\begin{Conjecture}
  Let $n$ and $c$ be fixed. Then there exists a graph $G$ on $O(r^{n-1})$ vertices such that $\chi(G)\geq n(c-1)+1$.
\end{Conjecture}

An affirmative answer to this question should reveal an interesting phenomenon. Consider the minimal number of vertices in a graph~$G$ such that $\chi(G)=N$ and $\lchi_r(G)\leq c$; then the rate of growth of this number of vertices (as the function in $r$) is $r^{\lfloor (N-1)/(c-1)\rfloor}$, that is --- it jumps at the values of $N$ congruent to $1$ modulo $c-1$.

In~\cite{kst}, the question is answered in affirmative for $n=2$ (using an example by Schmerl~\cite{schmerl}) and for $c=2$ (using an example related to Kneser graphs).
Another example verifying the question for $c=2$ was provided by Stiebitz~\cite{stie}. He  generalized and iterated the Mycielski construction obtaining, in particular, a series of graphs on $O(r^n)$ vertices with no odd cycles of length at most $2r+1$ such that $\chi(G)=n+2$. Stiebitz's proof utilizes a topological lemma by Lov\'asz~\cite{lovasz}. This proof is also reproduced in~\cite{gy-je-st}.

We also mention that Berlov and the author~\cite{berl-bogd} have obtained lower bounds for the number of vertices in a graph such that $\chi(G)\geq n$ and $\lchi_r(G)\leq 2$ for arbitrary values of $n$ and $r$. In a subsequent paper we will expand this estimate for arbitrary values of $c$.

The aim of this paper is to present an affirmative answer to the question above for all values of $n$ and $c$. Namely, we explicitly present a series of graphs verifying the following theorem.

\begin{Theorem}
  \label{th-lower}
  For every positive integers $c\geq 3$, $r$, and $n$ there exists a graph $G=(V,E)$ such that $\lchi_r(G)\leq c$, $\chi(G)\geq n(c-1)+1$, and
  \begin{equation}
    |V|=\frac{(2rc+1)^n-1}{2r}.
    \label{lower-my}
  \end{equation}
\end{Theorem}

The most difficult part of the justification is the proof of the lower bound for the chromatic number. This part is topological; it inspired by Stiebitz's proof mentioned above.

The necessary topological background is collected in Section~\ref{sub-top}. In Section~\ref{sub-constr} we describe the general construction; its properties are investigated in Section~\ref{sub-neighbor}. Finally, in section~\ref{proof} we prove Theorem~\ref{th-lower}.

\section{Topological background}
\label{sub-top}

Here we gather topological notions and facts needed in the sequel. We write $X\cong Y$ for homeomorphic topological spaces and $X\simeq Y$ for homotopy equivalent ones. For more detailed discussion see, e.g., \cite{kozlov}.

\subsection{Simplicial complexes}

An {\em (abstract) simplicial complex} is a pair $(V,{\mathsf K})$ where $V$ is a set and $\mathsf K\subseteq 2^V$ is a hereditary system of subsets of~$V$; this means that $F_1\subseteq F_2\in\mathsf K$ implies $F_1\in \mathsf K$. The elements of~$V$ are called {\em vertices}, and the sets in~$\mathsf K$ are called {\em simplices}. All simplicial complexes in our paper are finite, i.e. $|V|<\infty$. We will often denote a simplicial complex merely by~$\mathsf K$ assuming that its vertex set is $V(\mathsf K)=\bigcup \mathsf K$.

\smallskip
We say that $\phi\colon V(\mathsf K)\to \RR^d$ is a {\em geometric realization} of $\mathsf K$ if (\textit{i}) for every simplex $A=\{a_1,\dots,a_i\}\in \mathsf K$ the points $\phi(a_1),\dots,\phi(a_i)$ are affinely independent; and (\textit{ii}) for every two simplices $A,B\in\mathsf K$ we have $(\conv \phi(A))\cap (\conv \phi(B))=\conv \phi(A\cap B)$, where $\conv X$ is the convex hull of the set~$X$. It is known that each finite simplicial complex has a geometric realization. If $\phi$ is a geometric realization of $\mathsf K$, then we denote the topological subspace $\bigcup_{A\in\mathsf K}\conv \phi(A)\subset \RR^d$ by $\|\mathsf K\|$ and call it a {\em polyhedron} of $\mathsf K$. All polyhedra of~$\mathsf K$ are homeomorphic; thus this definition does not lead to an ambiguity.

Every point $\mathbf x\in\|\mathsf K\|$ lies in some simplex $\conv \phi(A)$ with $A\in\mathsf K$. The intersection of all such simplices in~$\mathsf K$ is again a simplex satisfying the same property. This simplex is called the {\em support} of $\mathbf x$ and denoted by~$\supp \mathbf x$; thus $x\in\conv \phi(\supp \mathbf x)$.

\smallskip
If $\mathsf K$ is a simplicial complex and $\mathsf L\subseteq \mathsf K$ is a hereditary subsystem of~$\mathsf K$ then we say that $\mathsf L$ is a {\em subcomplex} of $\mathsf K$ (we assume that the set of vertices of~$\mathsf L$ is $V(\mathsf L)=\bigcup \mathsf L$). If $\phi$ is a geometric realization of~$\mathsf K$, then $\phi\big|_{V(\mathsf L)}$ is also a geometric realization of~$\mathsf L$. In this case, we will always assume that the polyhedron $\|\mathsf L\|$ is a subspace of $\|\mathsf K\|$, i.e. $\|\mathsf L\|=\bigcup_{B\in\mathsf L}\phi(B)$. In particular, every simplex $A\in\mathsf K$ may be considered as the subcomplex $(A,2^A)$ of $\mathsf K$; thus we may write $\|A\|$ instead of $\conv\phi(A)$.

\subsection{Neighborhood complex and Lov\'asz's lemma}

The following important notion was introduced by Lov\'asz~\cite{lovasz}
.

\begin{Def}
  Let $G=(V,E)$ be a graph; we assume that it contains no isolated vertices. The {\em neighborhood complex} $\mathsf N(G)$ on the set of vertices $V$ consists of all subsets $A\subseteq V$ such that all elements of~$A$ have a common neighbor in~$G$ (this neighbor surely does not belong to~$A$).
\end{Def}

For instance, the neighborhood complex of the complete graph $K_r$ is an $(r-2)$-dimensional skeleton of the $(r-1)$-dimensional simplex; so $\|\mathsf N(K_r)\|\cong S^{r-2}$.

Lov\'asz has discovered a relation between the homotopy properties of~$\|\mathsf N(G)\|$ and the chromatic number of~$G$. To formulate this result, we need a notion of $k$-connectedness of a topological space.

As usual, we denote the unit ball and the unit sphere in~$\RR^d$ respectively by
$$
  B^d=\{{\mathbf x}\in\RR^d\colon |{\mathbf x}|\leq 1\}
  \qquad \text{and}\qquad
  S^{d-1}=\{\mathbf x\in\RR^d\colon |\mathbf x|=1\}.
$$
A nonempty topological space $X$ is {\em $k$-connected} if each continuous map $g\colon S^{m-1}\to X$ extends to a continuous map $\overline g\colon B^m\to X$, for $m=0,1,\dots,k+1$ (this condition for $m=0$ and $m=1$ means that $X$ is nonempty and path connected, respectively). It is well known that the sphere $S^k$ is $(k-1)$-connected but is not $k$-connected. Recall that homotopy equivalence preserves $k$-connectedness.

\begin{Lemma}[Lov\'asz, {\cite[Theorem 2]{lovasz}}]
  \label{lem-lovasz}
  Let $G$ be a graph. Assume that the polyhedron $\|{\mathsf N}(G)\|$ is $k$-connected. Then $\chi(G)\geq k+3$.
\end{Lemma}

This lemma was initially invented by Lov\'asz in order to find the chromatic number of Kneser graphs.

\subsection{Joins and nerves}

In view of Lov\'asz's lemma, it makes sense to seek for a graph $G$ such that $\|{\mathsf N}(G)\|$ is highly connected. For this, the following construction is useful.

The {\em join} of two topological spaces $X$ and $Y$ is defined as the quotient space
$$
  X*Y=(X\times Y\times [0,1])/\approx
$$
by the equivalence relation $\approx$ determined by $(x,y,0)\approx (x,y',0)$ and $(x,y,1)\approx (x',y,1)$ for all $x,x'\in X$ and $y,y'\in Y$. The join of $k$- and $\ell$-dimensional simplices is homeomorphic to a $(k+\ell+1)$-dimensional simplex. Moreover, $S^k*S^\ell\cong S^{k+\ell+1}$.

If $X$ and $Y$ are subspaces of~$\RR^m$ and $\RR^n$ then the join~$X*Y$ can be realized as a subspace of~$\RR^{m+n+1}$ in the following way. Choose two skew affine subspaces $U,V\subseteq \RR^{m+n+1}$ with $\dim U=m$, $\dim V=n$; we may regard $X$ and~$Y$ as the subspaces of~$U$ and~$V$, respectively. Then $X*Y\cong \{(1-t)x+ty\colon x\in X,\;y\in Y,\;t\in[0,1]\}$.

The {\em join} ${\mathsf K}*{\mathsf L}$ of two simplicial complexes $(U,\mathsf K)$ and $(V,\mathsf L)$ is defined as follows. We define $V(\mathsf K*\mathsf L)=U'\cup V'$, where $U'=U\times\{0\}$ and $V'=V\times\{1\}$ (thus we ensure that these sets are disjoint), and set
$$
  \mathsf K * \mathsf L=\bigl\{(A\times\{0\})\cup (B\times\{1\})\colon A\in{\mathsf K}, \; B\in{\mathsf L}\bigr\}.
$$

The two notions of a join agree in the sense that $\|\mathsf K*\mathsf L\|\cong \|\mathsf K\|*\|\mathsf L\|$; this is easily seen from the realization of the join described above.

\smallskip
In the sequel, for arbitrary graphs $G_1$ and $G_2$ we will construct a series of graphs $J_r$ such that $\|\mathsf N(J_r)\|\simeq \|\mathsf N(G_1)*\mathsf N(G_2)\|$. Let us introduce one more notion needed for the proof.

Let $X$ be a topological space, and let $\mathcal U=\{U(i)\colon i\in I\}$ be a covering of~$X$. The {\em nerve} of this covering is the simplicial complex with $I$ as the set of vertices; a subset $J\subseteq I$ is its simplex if $\bigcap_{i\in J} U(i)\neq\emp$. We will use the following well-known fact (see, e.g., \cite[Theorem 15.21]{kozlov}; we present only a particular case sufficient for our purposes).

\begin{Lemma}[Nerve lemma]
  \label{nerve}
  Let $\mathcal U=\{U_i\colon i\in I\}$ be a finite open covering of a compact metric space~$X$. Assume that for every $J\subseteq I$ the set $\bigcap_{i\in J} U(i)$ is either empty or contractible. Then the polyhedron of the nerve of~$\mathcal U$ is homotopy equivalent to~$X$.
\end{Lemma}


A nonempty set~$X\subseteq \RR^d$ is called {\em star-shaped} if there exists $a\in X$ such that for every $b\in X$ the whole segment $[a,b]$ lies in~$X$; in this case $a$ is called a {\em center} of~$X$. Obviously, each star-shaped set is contractible.

\section{The main construction}
\label{sub-constr}

Now we present the desired construction.

Let $G_1=(V_1,E_1)$ and $G_2=(V_2,E_2)$ be two graphs (we assume that they contain no isolated vertices), and let $r$ be a nonnegative integer. For $g_i\in V_i$, denote by $N_i(g_i)$ the set of all its neighbors in~$G_i$.

\begin{figure}[!th]
  \begin{center}
    \epsfig{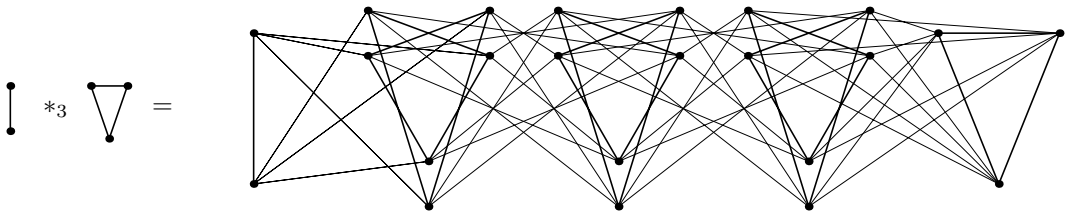}
    \caption{Graph $G_1*_r G_2$}
    \label{pic-join}
  \end{center}
\end{figure}

We define the graph $G_1 *_r G_2$ as follows. First, we define an auxiliary graph $J_r'=(U_r',F_r')$ by setting
\begin{gather*}
  U_r'=V_1\times V_2\times\{0,1,\dots,r+1\}, \\
  F_r'=\bigl\{\bigl((g_1,g_2,i),(g_1',g_2',j)\bigr)\in (U_r')^2\colon \;
    (g_1,g_1')\in E_1, \; (g_2,g_2')\in E_2, \;\text{and} \; |i-j|\leq 1\bigr\}.
\end{gather*}
The graph $J_r=G_1 *_r G_2=(U_r,E_r)$ is obtained by merging some vertices of the constructed graph. Namely, for every $g_1\in V_1$ we collapse all $|V_2|$ vertices of the form $(g_1,g_2,0)$ to a new vertex $(g_1,0)$, and for every $g_2\in V_2$ we collapse all $|V_1|$ vertices of the form $(g_1,g_2,r+1)$ to a new vertex $(g_2,r+1)$. Fig.~\ref{pic-join} shows a sample graph $K_2*_3 K_3$.

Notice that $G_1*_0G_2$ is just the usual join of graphs $G_1$ and $G_2$.

For every vertex $v\in U_r$, we denote by $N(v)$ the set of all neighbors of $v$ in $J_r$; these sets, together with all their subsets, form the complex $\mathsf N(J_r)$. For $r\geq 2$ these sets are
\begin{equation}
  \begin{aligned}
  N(g_1,0)&=\bigl(N_1(g_1)\times\{0\}\bigr)\cup \bigl(N_1(g_1)\times V_2\times\{1\}\bigr); \\
  N(g_1,g_2,1)&=\bigl(N_1(g_1)\times\{0\}\bigr)
    \cup\bigl(N_1(g_1)\times N_2(g_2)\times\{1,2\}\bigr); \\
  N(g_1,g_2,i)&=N_1(g_1)\times N_2(g_2)\times \{i-1,i,i+1\}\qquad \text{for $1<i<r$;} \\
  N(g_1,g_2,r)&=\bigl(N_2(g_2)\times\{r+1\}\bigr)
    \cup\bigl(N_1(g_1)\times N_2(g_2)\times\{r-1,r\}\bigr); \\
  N(g_2,r+1)&=\bigl(N_2(g_2)\times\{r+1\}\bigr)\cup \bigl(V_1\times N_2(g_2)\times \{r\}\bigr).
  \end{aligned}
  \label{N()}
\end{equation}

\section{The properties of the construction}
\label{sub-neighbor}

Now we will investigate the properties of the constructed graph $G_1*_r G_2$. First, we find the $r$-local chromatic number of $G_1 *_{2r} G_2$.

\begin{Lemma}
  \label{loc-join}
  For every graphs $G_1$ and $G_2$ and every positive integer $r$ we have
  $$
    \lchi_r(G_1*_{2r} G_2)=\max\{\lchi_r(G_1),\lchi_r(G_2)\}.
  $$
\end{Lemma}

\proof
  Since $G_1$ and $G_2$ are isomorphic to subgraphs of~$J_{2r}=G_1*_{2r} G_2$, we have $\lchi_r(J_{2r})\geq \max\{\lchi_r(G_1),\lchi_r(G_2)\}$.

  On the other hand, since the distance between~$V_1\times\{0\}$ and $V_2\times\{2r+1\}$ in~$J_{2r}$ is $2r+1$, each ball~$B$ of radius $r$ in $J_{2r}$ lies either in $J_{2r}\setminus(V_1\times\{0\})$ or in $J_{2r}\setminus(V_2\times\{2r+1\})$. Consider the first case. The projection to~$V_2$ is a graph homomorphism from $J_{2r}\setminus(V_1\times\{0\})$ to~$G_2$, and $\pr_{V_2}(B)$ is a ball of radius $r$ in~$G_2$; hence $\chi(B)\leq\chi(\pr_{V_2}(B))\leq \lchi_r(G_2)$. Similarly, in the second case we get $\chi(B)\leq \lchi_r(G_1)$, proving the converse inequality.
\qed

\medskip
Next, we deal with the neighborhood complex of $G_1 *_r G_2$.

\begin{Lemma}
  \label{join}
  For every graphs $G_1$ and $G_2$ and every integer $r\geq 2$ we have
  $$
    \|\mathsf N(G_1*_r G_2)\|\simeq \|\mathsf N(G_1)\|*\|\mathsf N(G_2)\|.
  $$
\end{Lemma}

\proof
Denote $\mathsf K=\mathsf N(G_1*_r G_2)$, $\mathsf M=\mathsf N(G_1)*\mathsf N(G_2)$. Let us construct a convenient geometric realization of~$\mathsf M$.
Consider some geometric realizations of $\mathsf N(G_1)$ and~$\mathsf N(G_2)$ in real spaces $R_1$ and~$R_2$; we may identify the vertices of $G_i$ with their images under these realizations. Now consider the space $R=R_1\times R_2\times \RR$; for every $g_1\in V_1$, identify the vertex $(g_1,0)\in \mathsf M$ with $(g_1,0,0)\in R$, and for every $g_2\in V_2$ identify the vertex $(g_2,1)\in\mathsf M$ with $(0,g_2,r+1)\in R$. This provides a geometric realization of~$\mathsf M$. For convenience, for every interval $I\subseteq \RR$ we denote by $R^I$ the ``strip'' $R_1\times R_2\times I\subseteq R$; in particular, $\|\mathsf M\|\subset R^{[0,\,r+1]}$.

Notice that for every nonempty topological subspaces $A_i\subseteq \|\mathsf N(G_i)\|$ the space $A_1*A_2$ may be regarded as a subset in $\|\mathsf M\|$. Moreover, for all $A_i,B_i\subseteq \|\mathsf N(G_i)\|$ such that $A_i\cap B_i\neq \emp$ we have $(A_1*A_2)\cap (B_1*B_2)=(A_1\cap B_1)*(A_2\cap B_2)$.

\smallskip
For every vertex $g_i\in V_i$ let us define the {\em star set} of vertex $g_i$ as
$$
  S_i(g_i)=\{x\in \|\mathsf N(G_i)\|\colon g_i\in\supp_{\mathsf N(G_i)} x\}\subset R_i.
$$
Each set $S_i(g_i)$ is open in~$\|\mathsf N(G_i)\|$. Next, for every $A_i\subseteq V_i$ the set $S_i(A_i):=\bigcap_{g_i\in A_i} S_i(g_i)$ is nonempty if and only if $A_i\in\mathsf N(G_i)$ (in fact, if $A_i\in \mathsf N(G_i)$ then $S_i(A_i)$ is the union of relative interiors of all simplices $\|B\|$ such that $A_i\subseteq B\in\mathsf N(G_i)$). For every $A_i\in \mathsf N(G_i)$, let us fix an arbitrary point $x_i(A_i)$ in the relative interior of $\|A_i\|$ (if $A_i=\{g_i\}$ then $x_i(A_i)=g_i$); one can easily see then that $S_i(A_i)$ is star-shaped with center~$x_i(A_i)$.

\begin{figure}[!th]
  \begin{center}
    \epsfig{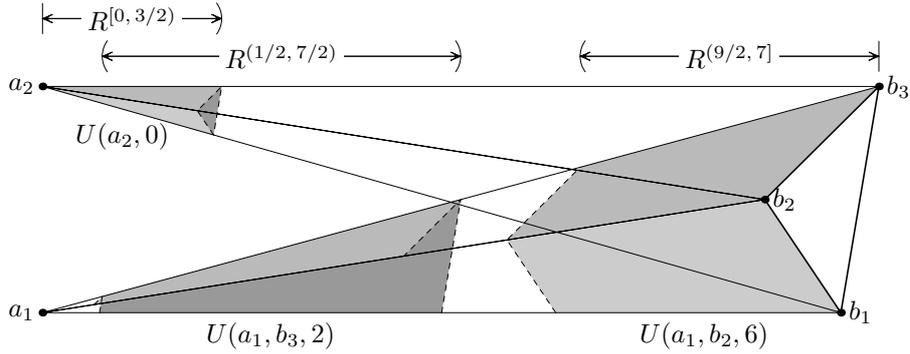}
    \caption{Covering $\mathcal U$ of $\|\mathsf N(G_1)*\mathsf N(G_2)\|$ (here $r=6$)}
    \label{pic-cover}
  \end{center}
\end{figure}

Now we construct a covering~$\mathcal U$ of~$\|\mathsf M\|$ with contractible intersections such that its nerve is~$\mathsf K$; by the Nerve lemma~\ref{nerve}, this implies the desired result. Set $\mathcal U=\{U(v)\colon v\in V(\mathsf K)\}$, where
\begin{gather*}
  U(g_1,0)=\bigl(S_1(g_1)*\|\mathsf N(G_2)\|\bigr)\cap R^{[0,\,3/2)}; \\
  U(g_2,r+1)=\bigl(\|\mathsf N(G_1)\|*S_2(g_2)\bigr)\cap R^{(r-1/2,\,r+1]}; \\
  U(g_1,g_2,i)=\bigl(S_1(g_1)*S_2(g_2)\bigr)\cap R^{(i-3/2,\,i+3/2)}
    \qquad \text{for $1\leq i\leq r$.}
\end{gather*}
Several sample sets $U(v)$ in $\|\mathsf N(K_2) * \mathsf N(K_3)\|$ are shown in Fig.~\ref{pic-cover}. Recall that in this example the underlying space~$R$ is 4-dimensional.

For an arbitrary maximal simplex $A\in \mathsf M$ we have $A=(A_1\times\{0\})\cup (A_2\times \{1\})$ with $A_i\in\mathsf N(G_i)$; then we have $\|A\|=\|A_1\|*\|A_2\|$, and one can easily see that this simplex is covered by the sets $U(g_1,g_2,t)$ with $g_i\in A_i$ and $t\in\{1,\dots,r\}$. Thus $\mathcal U$ is an open covering of~$\mathsf M$.

\smallskip
Let~$\mathsf L$ be the nerve of $\mathcal U$.
For every $A\subseteq V(\mathsf L)=V(\mathsf K)$ denote $U(A)=\bigcap_{v\in A} U(v)$. Now, one may verify that
\begin{gather*}
  U(A)\cap R^{\{0\}}=\begin{cases}
      S_1(\pr_{V_1}(A))\times\{0\}\times\{0\}, & \pr_{\RR}(A)\subseteq\{0,1\}; \\
      \emp, &\text{otherwise;}
    \end{cases} \\
  U(A)\cap R^{\{r+1\}}=\begin{cases}
      \{0\}\times S_2(\pr_{V_2}(A))\times\{r+1\}, & \pr_{\RR}(A)\subseteq\{r,r+1\}; \\
      \emp, &\text{otherwise;}
    \end{cases}\\
  U(A)\cap R^{(0,\,r+1)}=\begin{cases}
      \bigl(S_1(\pr_{V_1}(A))*\mathsf N(G_2)\bigr)\cap R^{(0,\,3/2)},
        &A\subseteq V_1\times\{0\}; \\
      \bigl(\mathsf N(G_1)*S_2(\pr_{V_2}(A))\bigr)\cap R^{(r-1/2,\,r+1)},
        &A\subseteq V_2\times\{r+1\}; \\
      \bigl(S_1(\pr_{V_1}(A))*S_2(\pr_{V_2}(A))\bigr)
        \cap R^{(a(A),\,b(A))},
        &\text{otherwise.}
    \end{cases}
\end{gather*}
Here we set $a(A)=\max\{0,\max(\pr_\RR (A))-3/2\}$ and $b(A)=\min\{r+1,\min(\pr_\RR (A))+3/2\}$. Next, the projection $\pr_{V_1}$ is defined as $\pr_{V_1}(g_1,0)=\pr_{V_1}(g_1,g_2,i)=g_1$ for $i=1,\dots,r$, and $\pr_{V_1}(g_2,r+1)=\emp$; the projection $\pr_{V_2}$ is defined similarly.

Now a straightforward check shows that $U(A)\cap R^{\{0\}}\neq \emp$ exactly if $A\subseteq N(g_1,0)$ for some $g_1\in V_1$, that $U(A)\cap R^{\{r+1\}}\neq \emp$ exactly if $A\subseteq N(g_2,r+1)$ for some $g_2\in V_2$, and that $U(A)\cap R^{(0,\,r+1)}\neq\emp$ exactly if $A\subseteq N(g_1,g_2,i)$ for some $g_1\in V_1$, $g_2\in V_2$, and $i\in\{1,\dots,r\}$. Thus $\mathsf L=\mathsf K$.

It remains to check that all nonempty sets of the form $U(A)$ are contractible; in fact, we will see that they are star-shaped. Assume that $U(A)\neq \emp$. If $U(A)\cap R^{\{0\}}\neq\emp$ then the set $U(A)$ is star-shaped with center $\bigl(x_1(\pr_{V_1}(A)),0,0\bigr)$. Similarly, if $U(A)\cap R^{\{r+1\}}\neq\emp$ then the set $U(A)$ is star-shaped with center $\bigl(0,x_2(\pr_{V_2}(A)),r+1\bigr)$. In the remaining case we have $U(A)=\bigl(S_1(\pr_{V_1}(A))*S_2(\pr_{V_2}(A))\bigr)\cap R^{(a(A),\,b(A))}$, and it is star-shaped with center $\bigl(x_1(\pr_{V_1}(A)),x_2(\pr_{V_2}(A)),(a(A)+b(A))/2\bigr)$.

\smallskip
Finally, applying the Nerve lemma~\ref{nerve} we get the required result.
\qed

\Remark. The statement of Lemma~\ref{join} remains valid for $r=1$ with essentially the same proof. On the other hand, for $r=0$ it does not hold in general; for instance, $\mathsf N(K_n)\cong S^{n-2}$, thus $\mathsf N(K_n *_0 K_m)\cong \mathsf N(K_{n+m})\cong S^{n+m-2}\not\simeq S^{n+m-3}\cong \mathsf N(K_n)*\mathsf N(K_m)$.

\section{Proof of Theorem~\ref{th-lower}}
\label{proof}

Let us fix the values of $c$ and $r$. We will use the induction on $n$ to construct the graph $G_n$ with $\bigl((2rc+1)^n-1\bigr)/(2r)$ vertices such that $\lchi_r(G_n)\leq c$ and $\mathsf N(G_n)\simeq S^{n(c-1)-1}$. Thus $\mathsf N(G_n)$ will be $(n(c-1)-2)$-connected; the required estimate then follows from Lemma~\ref{lem-lovasz}.

For $n=1$, the complete graph $K_c$ satisfies the desired properties since $|K_c|=c$ and $\mathsf N(K_c)\cong S^{c-2}$.

Now assume that $n>1$, and the graph $G_{n-1}$ is already constructed. Then we set $G_n=G_{n-1}*_{2r} K_r$. We have
$$
  |G_n|=(2rc+1)|G_{n-1}|+c=(2rc+1)\frac{(2rc+1)^{n-1}-1}{2r}+c=\frac{(2rc+1)^n-1}{2r}.
$$
Next, by Lemma~\ref{loc-join} we have $\lchi_r(G_n)=\max\{\lchi_r(G_{n-1}),\lchi_r(K_c)\}=c$. Finally, by Lemma~\ref{join}, we have $\mathsf N(G_n)\simeq \mathsf N(G_{n-1})*\mathsf N(K_c)
\simeq S^{(n-1)(c-1)-1}*S^{c-2}\cong S^{n(c-1)-1}$. The theorem is proved.
\qed

\bigskip
The author is grateful to S.L. Berlov who attracted author's attention to this question, and to R.N. Karasev and V.L. Dol'nikov for fruitful discussions.

\end{document}